\documentclass[11pt, a4paper]{article}

\usepackage{amsmath,amssymb,amsthm,amsfonts,color}
\usepackage{setspace}
\usepackage{mathrsfs}
\usepackage{graphicx}
\usepackage{tikz}

\newtheorem{theorem}{Theorem}

\newtheorem{lemma}{Lemma}

\usetikzlibrary{decorations.markings}
\tikzstyle{vertex}=[circle, draw, inner sep=2pt, minimum size=6pt]

\tikzstyle{filledvertex}=[circle, draw, fill, inner sep=2pt, minimum size=6pt]
\newcommand{\vertex}{\node[vertex]}

\tikzstyle{directed}=[postaction={decorate,
	decoration={markings,mark=at position 0.5 with {\arrow{stealth}}}
}]

\usepackage{geometry}
\begin{document}
\begin{spacing}{1.22}

\title{A note on hypergraph colorings\thanks{Supported by NSFC (No. 11601430) and China Postdoctoral Science Foundation (No. 2016M590969).}}

\author{\quad Yandong Bai\thanks{Corresponding author. E-mail address: bai@nwpu.edu.cn (Y. Bai).}\\[2mm]
\small Department of Applied Mathematics, Northwestern Polytechnical University, \\
\small Xi'an 710129, China}

\date{\today}
\maketitle

\begin{abstract}
Let $t\geqslant 2$ and $s\geqslant 1$ be two integers.
Define a $(t,s)$-coloring of a hypergraph
to be a coloring of its vertices using $t$ colors
such that each color appears on each edge at least $s$ times.
In this note,
we provide a sufficient condition for the existence of a $(t,s)$-coloring of a hypergraph by using the symmetric lopsided version of Lov\'asz Local Lemma.
Our result generalizes several known results on hypergraph colorings.

\smallskip
\end{abstract}

\section{Introduction}

A {\em hypergraph} $H$ consists of a collection $\mathcal{E}$ of subsets of a finite set $V$,
the members of $V$ are the {\em vertices}
and the members of $\mathcal{E}$ are the {\em edges} of $H$.
Define a {\em 2-coloring} of a hypergraph to be
a coloring of its vertices using two colors
such that no edge is monochromatic.

In 1975, by using the symmetric version of Lov\'{a}sz Local Lemma,
Erd\H{o}s and Lov\'{a}sz \cite{EL1975} obtained a sufficient condition for a hypergraph to have a 2-coloring.

\begin{theorem}[Erd\H{o}s and Lov\'{a}sz \cite{EL1975}]\label{Erdos-Lovasz theorem}
Let $H$ be a hypergraph in which every edge contains at least $k$ vertices
and meets at most $d$ other edges.
If
$$
e(d+1)\leqslant 2^{k-1},
$$
then $H$ has a 2-coloring.
\end{theorem}

By using the symmetric `lopsided' version of Lov\'{a}sz Local Lemma,
McDiarmid \cite{McDiarmid1997} improved Theorem \ref{Erdos-Lovasz theorem} in 1997 as follows.

\begin{theorem}[McDiarmid \cite{McDiarmid1997}]\label{McDiarmid theorem one}
Let $H$ be a hypergraph in which every edge contains at least $k$ vertices
and meets at most $d$ other edges.
If
$$
e(d+2)\leqslant 2^{k},
$$
then $H$ has a 2-coloring.
\end{theorem}

McDiarmid \cite{McDiarmid1997} also generalized $2$-colorings to {\em $t$-colorings},
that is, a coloring of its vertices using $t$ colors such that each color appears on each edge.
An analogous sufficient condition has been obtained.

\begin{theorem}[McDiarmid \cite{McDiarmid1997}]\label{McDiarmid theorem two}
Let $H$ be a hypergraph in which every edge contains at least $k$ vertices
and meets at most $d$ other edges.
Let $t\geqslant 2$ be an integer.
If
$$
e\left(1-\frac{1}{t}\right)^{k}\Big((d+1)(t-1)+1\Big)\leqslant 1,
$$
then $H$ has a $t$-coloring.
\end{theorem}

Motivated by the above results,
we consider a generalization of $t$-colorings.
Define a {\em $(t,s)$-coloring} of a hypergraph
to be a coloring of its vertices using $t$ colors
such that each color appears on each edge at least $s$ times.
Note that the $(t,1)$-coloring is the general $t$-coloring.
In 2010,
Chen et al. \cite{CDM2010} considered $(t,2)$-colorings of a hypergraph
and showed the following result.

\begin{theorem}[Chen et al. \cite{CDM2010}]\label{Chen et al. theorem}
Let $H$ be a hypergraph in which every edge contains at least $k$ vertices
and meets at most $d$ other edges.
Let $t\geqslant 2$ be an integer with $k\geqslant 2t$.
If
$$
e\left(1-\frac{1}{t}\right)^{k-1}\left(1-\frac{1}{t}+\frac{k}{t}\right)\Big((d+1)(t-1)+1\Big)\leqslant 1
$$
then $H$ has a $(t,2)$-coloring.
\end{theorem}

In this note,
by using the symmetric lopsided version of Lov\'asz Local Lemma,
we provide a sufficient condition for the existence of a $(t,s)$-coloring of a hypergraph for general $t\geqslant 2$ and $s\geqslant 1$.

\begin{theorem}\label{main theorem}
Let $H$ be a hypergraph in which every edge contains at least $k$ vertices
and meets at most $d$ other edges.
Let $t\geqslant 2$ and $s\geqslant 1$ be two integers with $k\geqslant st$.
If
\begin{equation*}
e\sum_{j=0}^{s-1}\binom{k}{j}\left(\frac{1}{t}\right)^{j}\left(1-\frac{1}{t}\right)^{k-j}\Big((d+1)(t-1)+1\Big)\leqslant 1,
\end{equation*}
then $H$ has a $(t,s)$-coloring.
\end{theorem}

\noindent
\textbf{Remark 1.}
Taking $t=2$ and $s=1$ in Theorem \ref{main theorem},
we get Theorem \ref{McDiarmid theorem one};
taking $s=1$ in Theorem \ref{main theorem},
we get Theorem \ref{McDiarmid theorem two};
and taking $s=2$ in Theorem \ref{main theorem},
we get Theorem \ref{Chen et al. theorem}.
For more results on hypergraph colorings and related problems,
we refer the reader to \cite{Bretto2013,FM2011,FM2013,Halldorsson2010,HY2013,HY2014,LRS2011,NI2008,Seymour1974,Szymanska2012,Vishwanathan2003}.

\section{Proof of Theorem \ref{main theorem}}

The famous Lov\'asz Local Lemma will play a key role in our proof.

\begin{lemma}[Local Lemma: symmetric lopsided version, see in \cite{AS2008}]\label{local lemma}
Let $A_{1},\ldots,A_{n}$ be events in a probability space $\Omega$ with a lopsidependency graph $G$.
Suppose that each event $A_{i}$ has probability at most $p$
and each vertex in $G$ has degree at most $\Delta$.
If
\begin{equation*}
ep(\Delta+1)\leqslant 1
\end{equation*}
then $Pr[\bigwedge_{i=1}^{n}\overline{A}_{i}]>0$.
\end{lemma}

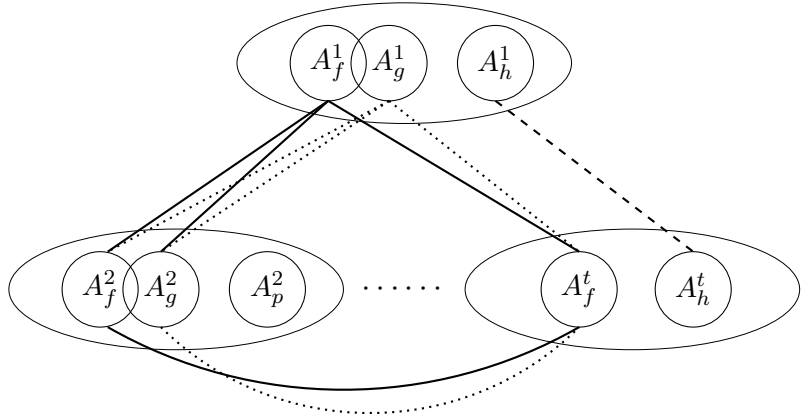
\begin{figure}[ht]
\begin{center}
\begin{tikzpicture}
\tikzstyle{vertex}=[circle,inner sep=2pt, minimum size=0.1pt]

\draw (-3,-1.5) [line width=0.2pt] ellipse (2.2 and 0.8);
\draw (0,1.5) [line width=0.2pt] ellipse (2.2 and 0.8);
\draw (3,-1.5) [line width=0.2pt] ellipse (2.2 and 0.8);
\draw (-1,1.5) circle (0.5);
\draw (-0.2,1.5) circle (0.5);
\draw (1.2,1.5) circle (0.5);

\vertex  at (-1,1.5)[label=center:$A_{f}^{1}$]{};
\vertex  at (-0.2,1.5)[label=center:$A_{g}^{1}$]{};
\vertex  at (1.2,1.5)[label=center:$A_{h}^{1}$]{};
\draw (-4,-1.5) circle (0.5);
\draw (-3.2,-1.5) circle (0.5);
\draw (-1.8,-1.5) circle (0.5);

\vertex  at (-4,-1.5)[label=center:$A_{f}^{2}$]{};
\vertex  at (-3.2,-1.5)[label=center:$A_{g}^{2}$]{};
\vertex  at (-1.8,-1.5)[label=center:$A_{p}^{2}$]{};
\draw (2.3,-1.5) circle (0.5);
\draw (3.8,-1.5) circle (0.5);

\vertex  at (2.3,-1.5)[label=center:$A_{f}^{t}$]{};
\vertex  at (3.8,-1.5)[label=center:$A_{h}^{t}$]{};

\vertex  at (0,-1.5)[label=center:{\large $\cdots \cdots$}]{};
\draw [line width=0.8pt](-1,1)--(-3.9,-1);
\draw [line width=0.8pt](-1,1)--(-3.2,-1);
\draw [line width=0.8pt](-1,1)--(2.3,-1);

\draw [dotted,line width=0.8pt] (-0.2,1)--(-3.9,-1.0);
\draw [dotted,line width=0.8pt](-0.2,1)--(-3.2,-1);
\draw [dotted,line width=0.8pt](-0.2,1)--(2.3,-1);

\draw [dashed,line width=0.8pt](1.2,1)--(3.8,-1);


\vertex (a) at (-3.9,-2)[]{};
\vertex (b) at (-3.2,-2)[]{};
\vertex (c) at (-1.8,-2)[]{};

\draw [line width=0.8pt] (a) arc (240: 300: 6.2);

\draw [dotted,line width=0.8pt] (b) arc (225: 315: 3.89);


\end{tikzpicture}
\end{center}
\caption{An illustration of the constructed graph $G$.}
\label{figure: lopsidependency graph}
\end{figure}

Color the vertices of $H$ by using colors $1,\ldots,t$ independently and with equal probability.
Let $A_{f}^{i}$ be the event that
the color $i$ appears on the edge $f$ at most $s-1$ times.
Then
\begin{equation*}
Pr[A_{f}^{i}]
= \sum_{j=0}^{s-1}\binom{|f|}{j}\left(\frac{1}{t}\right)^{j}\left(1-\frac{1}{t}\right)^{|f|-j}.
\end{equation*}
Let 
$$
\theta(|f|,j)=\binom{|f|}{j}\left(\frac{1}{t}\right)^{j}\left(1-\frac{1}{t}\right)^{|f|-j}.
$$
By $t\geqslant 2$, we have $\theta(|f|,j)>0$.
Since $|f|\geqslant k\geqslant st> jt$ for each $j\in \{0,\ldots,s-1\}$,
we have
\begin{equation*}
\begin{split}
\frac{\theta(|f|,j)}{\theta(|f|+1,j)}
&=\frac{\binom{|f|}{j}\left(\frac{1}{t}\right)^{j}\left(1-\frac{1}{t}\right)^{|f|-j}}
{\binom{|f|+1}{j}\left(\frac{1}{t}\right)^{j}\left(1-\frac{1}{t}\right)^{|f|+1-j}}\\
&=\frac{1-\frac{j}{|f|+1}}{1-\frac{1}{t}}\\
&> 1.
\end{split}
\end{equation*}
Hence $\theta(|f|,j)$ is a decreasing function and, by $|f|\geqslant k$,
we have $\theta(|f|,j)\leqslant \theta(k,j)$ for each edge $f$ and each $j\in \{0,\ldots,s-1\}$.
It follows that
\begin{equation}
\begin{split}
Pr[A_{f}^{i}]
&= \sum_{j=0}^{s-1}\theta(|f|,j)\leqslant \sum_{j=0}^{s-1}\theta(k,j)=\sum_{j=0}^{s-1}\binom{k}{j}\left(\frac{1}{t}\right)^{j}\left(1-\frac{1}{t}\right)^{k-j}.
\end{split}
\end{equation}
We construct a $t$-partite graph $G$ with 
$$
V(G)=\{A_{f}^{i}:~i\in \{1,\ldots,t\},~f\in E(H)\},~~E(G)=\{A_{f}^{i}A_{g}^{i}:~i\neq j,~f\cap g\neq \emptyset\},
$$ 
see an illustration in Figure \ref{figure: lopsidependency graph}.
Since each edge meets at most $d$ other edges and there are $t$ colors in total,
we have
\begin{equation}
\Delta(G)\leqslant (d+1)(t-1).
\end{equation}
Now we show that $G$ is a lopsidependency graph for the events $A_{f}^{j}$.
The following powerful lemma is needed.

\begin{lemma}[McDiarmid \cite{McDiarmid1997}]\label{lemma: McDiarmid}
Let $n,t$ be two positive integers
and let $P$ be a product measure on $\Omega=\{1,\ldots,n\}^{t}$.
For each (color) $j\in \{1,\ldots,t\}$,
let $\mathscr{C}^{(j)}$ be a collection of subsets of $V=\{1,\ldots,n\}$.
For each set $f\in \mathscr{C}^{(j)}$,
let $D_{f}^{(j)}$ be a hereditary collection (downset) of subsets of $f$
and let $A_{f}^{(j)}$ be the event
$$
A_{f}^{(j)}=\{\omega\in \Omega:~\{i\in f: \omega_{i}=j\}\in D_{f}^{j}\}.
$$
Let $\widehat{G}$ be the $t$-partite graph with parts $\{A_{f}^{(j)}\}$ for each $j\in \{1,\ldots,t\}$
and satisfying that vertices $A_{f}^{(j)}$ and $A_{f'}^{(j')}$ are adjacent if and only if
$j\neq j'$ and $f\cap f'\neq \emptyset$.
Then $\widehat{G}$ is a lopsidependency graph for the events $A_{f}^{(j)}$.
\end{lemma}

In our $(t,s)$-coloring,
for each $j\in \{1,\ldots,t\}$,
let $D_{f}^{(j)}$ be the set of all subsets of $f$ of size at most $s-1$.
Then one can see that the event $A_{f}^{(j)}$ defined in Lemma \ref{lemma: McDiarmid} is the same as the event $A_{f}^{j}$ defined above.
So by Lemma \ref{lemma: McDiarmid} the graph $G$ is a lopsidependency graph for the events $A_{f}^{j}$.
The required result therefore follows directly from Inequalities (1), (2) and Lemma \ref{local lemma}.

\end{spacing}

\begin{thebibliography}{0}

\bibitem{AS2008}
N. Alon, J. H. Spencer,
The Probabilistic Method, 3rd Edition,
Wiley, 2008.

\bibitem{Bretto2013}
A. Bretto,
Hypergraph theory,
Springer International Publishing, 2013.

\bibitem{CDM2010}
X. Chen, Z. Du, J. Meng,
Coloring and the Lov\'{a}sz Local Lemma,
Applied Math. Letters 23 (2010) 219-221.

\bibitem{EL1975}
P. Erd\H{o}s, L. Lov\'{a}sz,
Problems and results on 3-chromatic hypergraphs and some related questions,
in: A. Hajnal et al., eds., Infinite and Finite Sets (North-Halland, Amsterdam, 1975) 609-628.

\bibitem{FM2011}
A. Frieze, P. Melsted,
Randomly coloring simple hypergraphs,
Information Processing Letters 111 (2011) 848-853.

\bibitem{FM2013}
A. Frieze, D. Mubayi,
Coloring simple hypergraphs,
J. Combin. Theory Ser. B 103 (2013) 767-794.

\bibitem{Halldorsson2010}
M. M. Halld\'{o}rsson,
Online coloring of hypergraphs,
Information Processing Letters 110 (2010) 370-372.

\bibitem{HY2013}
M. A. Henning, A. Yeo,
2-colorings in $k$-regular $k$-uniform hypergraphs,
European J. Combin. 34 (2013) 1192-1202.

\bibitem{HY2014}
M. A. Henning, A. Yeo,
On 2-colorings of hypergraphs,
J. Graph Theory 80 (2015) 112-135.

\bibitem{LRS2011}
X. Li, A. Rudra, R. Swanminathan,
Flexible coloring,
Information Processing Letters 111 (2011) 538-540.

\bibitem{McDiarmid1997}
C. McDiarmid,
Hypergraph colouring and the Lov\'{a}sz Local Lemma,
Discrete Math. 167/168 (1997) 481-486.

\bibitem{NI2008}
J. Nagy-Gy\"{o}rgy, Cs. Imreh,
Online hypergraph coloring,
Information Processing Letters 109 (2008) 23-26.

\bibitem{Seymour1974}
P. D. Seymour,
On the two coloring of hypergraphs,
Quart. J. Math. Oxford Ser. 25 (1974) 303-312.

\bibitem{Szymanska2012}
E. Szyma\'{n}ska,
$H$-colorings of dense hypergraphs,
Information Processing Letters 112 (2012) 899-902.

\bibitem{Vishwanathan2003}
S. Vishwanathan,
On 2-coloring certain $k$-uniform hypergraphs,
J. Combin. Theory Ser. A 101 (2003) 168-172.

\end{thebibliography}
\end{document}